\documentclass{amsart}

\usepackage{amsmath}

\theoremstyle{plain}
\newtheorem{thm}{Theorem}[section]
\newtheorem{lem}[thm]{Lemma}

\theoremstyle{definition}
\newtheorem{defn}{Definition}[section]
\newtheorem{exmpl}[defn]{Example}

\newcommand{\ul}{{\mathcal L}}
\newcommand{\fdf}[1]{{\emph{ #1}}}
\newcommand{\conf}{{\mathcal C}}
\newcommand{\sconf}{{\mathcal S}{\mathcal C}}
\newcommand{\vrtcs}{{\mathcal V}}
\newcommand{\edge}{{\mathcal E}}
\newcommand{\KM}{Kapovich and Millson}
\newcommand{\mand}{\text{ and }}
\newcommand{\mif}{\text{ if } }

\newcommand{\qf}{quasifunctional}
\newcommand{\qfl}{quasifunctional linkage}
\newcommand{\compose}{\circ}
\newcommand{\PC}{Peaucellier}
\newcommand{\capshun}[2]{\caption{ #2}\label{fig:#1}}
\newcommand{\commentout}[1]{}

\DeclareMathOperator{\Cl}{Cl}   
\newcommand{\Eu}{\mathrm {Euc}}
\newcommand{\Euc}{\mathrm {Euc}}
\newcommand{\Tran}{\mathrm {Tran}}
\newcommand{\lineseg}[2]{\overline{#1#2}}
\newcommand{\thma}{Theorem \ref{thm:a}}
\newcommand{\thmb}{Theorem \ref{thm:b}}
\newcommand{\thmc}{Theorem \ref{thm:c}}

       \input boxedeps.tex 
       \input boxedeps.cfg 

\begin{document}

\title{Semiconfiguration Spaces of Planar Linkages}
\author[KING]{Henry C. King}
\address{
Department of Mathematics\\
University of Maryland\\
College Park, Maryland, 20742}
\email{hking@math.umd.edu}

\begin{abstract}
This paper characterizes which subsets of ${\mathbb C}^n$ can be 
the set of positions of $n$ points on a linkage in ${\mathbb C}$.
For example, assuming compactness they are just 
compact semialgebraic sets.
Noncompact configuration spaces are semialgebraics sets
invariant under the Euclidean group, with compact quotient.
\end{abstract}

\maketitle

\section{Linkages}

Loosely speaking, a linkage is an ideal mechanical device
consisting of a bunch of stiff rods sometimes attached at their ends
by rotating joints. A planar realization is some way of
placing this linkage in the plane.
The configuration space for a linkage is the space of all such planar realizations,
which can be determined by looking at all possible positions
of the ends of the rods.
Such configurations spaces have been studied for example in
\cite{KM} and \cite{K}.
In this paper, we look at semiconfiguration spaces of linkages,
where we look at all possible positions of only some of the
points on the linkage.
  For example, what curve does a particular point on the linkage
  trace out?
   We give a complete description of possible semiconfiguration spaces
in \thma\ below.
For example, compact semiconfiguration spaces correspond exactly
to compact semialgebraic sets.

Suppose $L$ is a finite one dimensional simplicial
complex, in other words, a finite set $\vrtcs(L)$
of vertices and a finite set $\edge(L)$ of edges
between certain pairs of vertices.
An \fdf{abstract linkage} is a finite one dimensional simplicial
complex $L$ with a positive number
$\ell ({\lineseg vw})$ assigned to each edge ${\lineseg vw}$,
i.e., a function $\ell \colon \edge(L)\to (0,\infty )$.
A \fdf{planar realization} of an abstract linkage $(L,\ell )$ is a mapping
$\varphi \colon \vrtcs(L)\to {\mathbb C}$ so that 
$|\varphi (v)-\varphi (w)|=\ell ({\lineseg vw})$ 
for all  edges ${\lineseg vw}$.

We will often wish to fix some of the vertices of a linkage whenever we
take a planar realization.
So we say that a \fdf{planar linkage} $\ul$ is a foursome $(L,\ell ,V,\mu )$ where
 $(L,\ell )$ is an abstract linkage, $V\subset \vrtcs(L)$ is a subset of its vertices, and $\mu \colon V\to {\mathbb C}$.
 So $V$ is the set of fixed vertices and $\mu $ tells where to fix them.
The configuration space of realizations is defined by:
$$\conf(\ul)=\left\{\,\varphi \colon \vrtcs(L)\to {\mathbb C}\Biggm|
\begin{array}{cll}
\varphi (v)=\mu (v) & \mif\ v\in V \\
   |\varphi (v)-\varphi (w)|
   =\ell ({\lineseg vw})&\text{ for all edges } 
    {\lineseg vw}\in \edge(\ul)
    \end{array}
\,\right\}$$   
If $W\subset \vrtcs(L)$ is a collection of vertices, 
the semiconfiguration space is the set of restrictions of realizations to $W$,
$$\sconf(\ul,W)=\{\,\varphi \colon W\to {\mathbb C} \mid 
\text{there is a } \varphi '\in \conf(\ul) \text{ so that }
 \varphi =\varphi '|_W  \,\}$$

If we order the elements of $W$ as $w_1,w_2,\ldots ,w_k$ then there
is a natural identification of $\sconf(\ul,W)$ with a subset of ${\mathbb C}^k$
where $\varphi $ is identified with the point $(\varphi (w_1),\ldots ,\varphi (w_k))$.
With this identification, we see that the semiconfiguration space of
a linkage is the projection of its configuration space to some coordinate
subspace.

Note that $\conf(\ul)$ is a real algebraic set inside ${\mathbb C}^n$ since it is
given by polynomial equations of the form $z_i=$ a constant
and $|z_i-z_j|^2=$ a constant.
So $\sconf(\ul,W)$ is the projection of a real algebraic set.
By the Tarski-Seidenberg theorem \cite{S}, 
projections of real algebraic sets are semialgebraic sets.
A \fdf{semialgebraic set} is a finite union of sets of the form
$$\{\,x\in {\mathbb R}^n\mid p_i(x)=0, q_j(x)\ge 0, \mand r_k(x)>0\,\}$$
for collections of polynomials $p_i$, $q_j$ and $r_k$.
In other words it is the closure 
under the boolean operations of finite union, intersection,
and complement,
of the family sets of the form 
$p^{-1}([0,\infty ))$ , where $p$ is a polynomial.


We are now ready to completely characterize semiconfiguration spaces.
Up to linear maps, they are just compact semialgebraic sets
cartesian product with ${\mathbb C}^m$.
But we can be even more precise.

Let $\Eu(2)$ denote the group of Euclidean motions of ${\mathbb C}$.
So a general element of $\Eu(2)$ is a map
$z\mapsto \omega z+z_0$ or $\omega \overline z+z_0$ where
$\omega \in {\mathbb C}$ satisfies $|\omega |=1$.
We say a subset $Z\subset {\mathbb C}^k$ is \fdf{virtually compact}
if $Z$ is either compact, or it is invariant under the diagonal action
of $\Eu(2)$, with compact quotient.

\begin{thm}\label{thm:a}
Suppose $X\subset {\mathbb C}^n$.
\begin{enumerate}
\item The following are equivalent:
\begin{enumerate}
\item There is a linkage $\ul$ and a $W\subset \vrtcs(\ul)$
so that $\sconf(\ul,W)=X$.
\item After perhaps permuting the coordinates,
 $X=Y_1\times  Y_2 \times \ldots \times Y_m $
where each 
$Y_i\subset {\mathbb C}^{k_i}$ is  a virtually compact semialgebraic set.
\end{enumerate}
\item Moreover, for connected linkages the following are equivalent:
\begin{enumerate}
\item There is a connected linkage $\ul$ and a $W\subset \vrtcs(\ul)$
so that $\sconf(\ul,W)=X$.
\item  $X$ is a virtually compact semialgebraic set.
\end{enumerate}
\end{enumerate}
\noindent As an extra bit of information, if the connected linkage
$\ul$ has any fixed vertices or if $W$ is empty, then $\sconf(\ul,W)$
is compact.
Otherwise, $\sconf(\ul,W)$ is invariant under the diagonal action
of $\Eu(2)$, with compact quotient, but $\sconf(\ul,W)$ itself
is noncompact.
In the general case if $\ul$ is perhaps disconnected,
then $\sconf(\ul,W)$ is compact if and only if each
component of $\ul$ which contains points of $W$
also contains a fixed vertex.
\end{thm}

\begin{exmpl}\label{ex:1}
We will illustrate \thma\ with an example.
Suppose that $X\subset {\mathbb C}^2$ is the diagonal,
$X=\{\,(z_1,z_2)\in {\mathbb C}^2 \mid z_1=z_2\,\}$.
Then $X$ is $\Euc(2)$ invariant and $X/\Euc(2)$ is a single point.
Let $\ul$ be the linkage with five vertices $A,B,C,D,E$
and with edges $\lineseg AB,\lineseg BC,\lineseg AC$ of length 1 and edges 
$\lineseg AD,\lineseg BD,\lineseg CD,\lineseg AE,\lineseg BE,\lineseg CE$ 
of length $\sqrt 5/4$.
Let $W=\{D,E\}$.
We do not fix any vertices of $\ul$.
Then in any realization of $\ul$, the vertices $A,B,C$
form an equilateral triangle and the vertices $D$ and $E$
both must lie at the barycenter.
Consequently, $\sconf(\ul,W)=X$.
Note that the configuration space $\conf(\ul)$
is higher dimensional, since the triangle $ABC$ can
rotate around $D$ and $E$.
In fact $\conf(\ul)$ is a single orbit
of $\Euc(2)$ (with no isotropy).
\end{exmpl}

We can further characterize semiconfiguration spaces
according to the number of fixed vertices.
For clarity, we restrict attention to the connected case.
The corresponding statement for disconnected linkages
follows from the observation that the
semiconfiguration space of the disjoint union of two
linkages is the cartesian product of their
semiconfiguration spaces, (see Lemma \ref{lem:6}).

\begin{thm}\label{thm:b}
Let $Z\subset {\mathbb C}^n$ be a virtually compact semialgebraic set,
and suppose that $Z$ is invariant under the action of a
subgroup $G$ of $\Eu(2)$, with compact quotient.
Then there is a connected linkage $\ul$ and a $W\subset \vrtcs(\ul)$ so that
$\sconf(\ul,W)=Z$ and so that:
\begin{enumerate}
\item If $G$ is $\Eu(2)$,
then  $\ul$ has no fixed vertices.
\item If $G$ is conjugate to $O(2)$,
then  $\ul$ has only one fixed vertex,
and that vertex is fixed at the fixed point of $G$.
\item If $G$ is an order two subgroup generated
by a reflection,
then  $\ul$ has only two fixed vertices,
and these vertices are fixed at points on the fixed line of $G$.
\item Otherwise,  $\ul$ has only three fixed vertices.
\end{enumerate}
\noindent
Moreover, we have a converse to this result.
Suppose $\ul$ is a connected linkage with $m$ fixed vertices
and $W\subset \vrtcs(\ul)$,
then:
\begin{enumerate}
\item[5.] If $m=0$, then $\sconf(\ul,W)$ is a semialgebraic set
invariant under
the action of $\Eu(2)$, with compact quotient.
\item[6.] If $m=1$, then $\sconf(\ul,W)$ is a compact 
semialgebraic set invariant under
the action of a subgroup conjugate
to $O(2)$.
\item[7.] If $m=2$, then $\sconf(\ul,W)$ is a compact 
semialgebraic set invariant under
the action of the order two subgroup generated
by a reflection.
\item[8.] If $m>0$, then $\sconf(\ul,W)$ is a compact
semialgebraic set.
\end{enumerate}
\end{thm}


\section{Constructing polynomial \qf\ linkages}

A linkage $\ul$ is \fdf{\qf\ } for a map $f\colon {\mathbb C}^n\to {\mathbb C}^m$ 
 if there are 
 vertices $w_1,\ldots ,w_n$ and $v_1,\ldots ,v_m$ of $\ul$ so that
 if $p\colon \conf(\ul)\to {\mathbb C}^m$ 
 and $q\colon \conf(\ul)\to {\mathbb C}^n$ are the maps
 $p(\varphi )=(\varphi (v_1),\ldots ,\varphi (v_m))$ and
$q(\varphi )=(\varphi (w_1),\ldots ,\varphi (w_n))$, 
then $p=f\compose q$.
We call $w_1,\ldots ,w_n$ the input vertices and
$v_1,\ldots ,v_m$ the output vertices.
So in other words $\ul$ is \qf\ if
the output vertices are $f$ of the input vertices.
The \fdf{domain} of a \qfl\ is $q(\conf(\ul))$.
In general, repetitions of input and output vertices are allowed.
But for convenience, for all quasifunctional linkages in this paper,
we will assume there
are no repetitions, i.e., $v_i\neq v_j$ and $w_i\neq w_j$ if $i\neq j$.

%

In \cite{KM} or \cite{K}, \qf\ linkages were constructed
for any real polynomial map ${\mathbb C}^n\to {\mathbb C}^m$.
In fact these linkages had some stronger properties,
which we do not need in this paper.
We will reproduce these constructions here,
 simplified when appropriate.
Essentially, polynomial \qf\ linkages were constructed
in the nineteenth century, although \KM\
pointed out some necessary corrections to the old constructions.
In particular, whenever there is a rectangle in a linkage,
it should be rigidified by adding another edge joining the
midpoints of two opposite edges.
This will prevent certain degenerate realizations which
destroy \qf ity.
In the diagrams below, we represent this rigidifying edge by 
a gray line.
The second correction made in \cite{KM} is in the \PC\ inversor below,
which we correct here by adding a simulated cable (see below).
One could also add a simulated telescoping edge as was done in 
\cite{KM}.

\begin{thm}[\KM]\label{thm:c}
For any real polynomial function $f\colon {\mathbb C}^n\to {\mathbb C}^m$
and any compact $K\subset {\mathbb C}^n$ there is a \qfl\ $\ul$
for $f$ whose domain contains $K$,
and whose input and output vertices are all distinct.
\end{thm}

The proof of \thmc\ will occupy the rest this section.
We now make some observations which reduce the proof of
\thmc\ to some special cases.
\begin{itemize}
\item
The first observation is that by taking two such
\qf\ linkages and  attaching the inputs of one to the outputs
of the other, we obtain a \qf\ linkage for the composition.
Consequently, it suffices to find \qf\ linkages 
for the elementary operations of addition, multiplication,
and complex conjugation.
Since $zw=(z+w)^2/4-(z-w)^2/4$ we may replace multiplication
by squaring and real scalar multiplication.
\item
The next observation is that there is a nontrivial \qfl\
for the identity with  domain all of ${\mathbb C}$.
Recall Example \ref{ex:1} with semiconfiguration space
the diagonal in ${\mathbb C}^2$.
With $D$ as input and $E$ as output, this is a
\qfl\ for the identity with distinct input and
output vertices.  Consequently, by attaching these to the
input and output vertices of any \qfl, we may transform any
\qfl\ to one with distinct input and output vertices.
\end{itemize}

\subsection{Simulating interior joints, cables, and telescoping edges}

In our model of linkages, edges are connected only at their ends.
Actual linkages used in real life might have a connection in the
middle of an edge.
This may be simulated as in Figure \ref{fig:1}.
In any realization, $C$ must lie on the line segment from $A$ to $B$.
Thus when drawing linkages, it is allowable to draw
a joint in the middle of an edge.

\begin{figure}
\BoxedEPSF {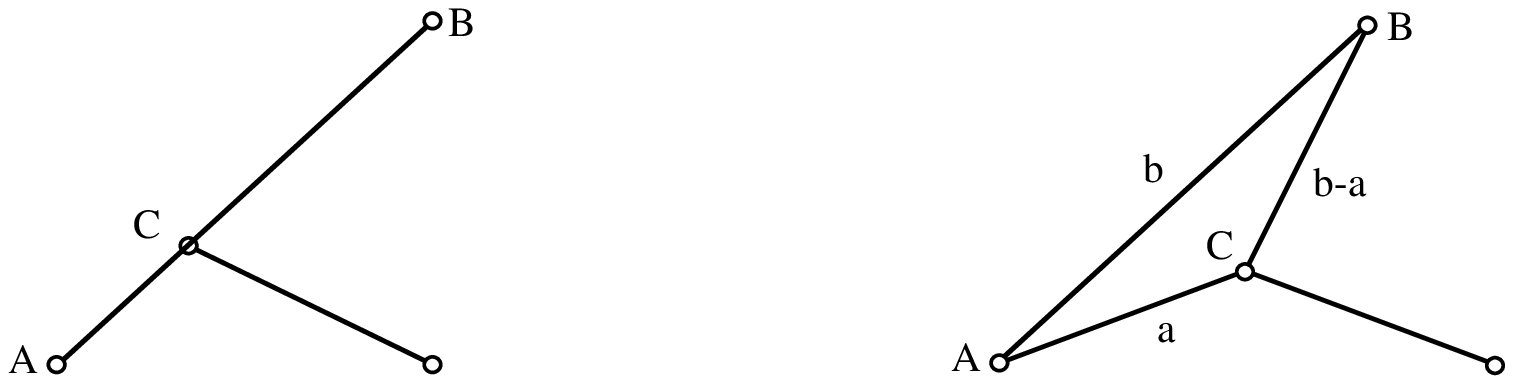 scaled 700}
\capshun 1 { How to put a joint in the middle of an edge}
\end{figure}

Although we will not use the following constructions in this paper
in an essential way,
we point out that using semiconfiguration spaces, we can also
simulate other types of linkages.
For example, suppose we want two vertices $A$ and $B$ connected by a cable, so 
the distance between them is constrained to be $\le b$.
More generally, suppose we wish to connect $A$ and $B$
by a telescoping edge,
so the distance between them is constrained to be
in the interval $[a,b]$.
This can be simulated as in Figure \ref{fig:2}.  Since we are using semiconfiguration
spaces, we can ignore the position of the vertex $D$.
To simulate a cable, we take $c=d=b/2$.
To simulate a telescoping edge with $0<a<b$, we take 
$c=(a+b)/2$, $d=(b-a)/2$.

\begin{figure}
\BoxedEPSF {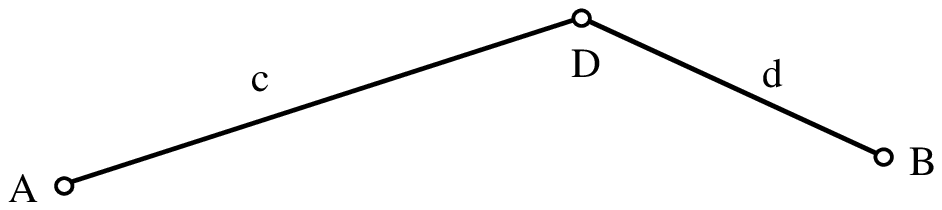 scaled 700}
\capshun 2 { Simulating a cable or telescoping edge}
\end{figure}

\subsection{The useful pantograph}

We have reduced \thmc\ to finding
 \qf\ linkages for $(z,w)\mapsto (z+w)/2$, $z\mapsto \lambda z$
for $\lambda $ real,
$z\mapsto z^2$, and $z\mapsto \overline z$,
all with domain containing an arbitrarily large compact set.

The first two functions can all be obtained from one type of linkage,
the pantograph shown in Figure \ref{fig:3}.
It is a rigidified rectangle $DEBF$ with two extended sides.
Because of the rigidification, in any realization the line
$\lineseg AD$ is parallel to $\lineseg BF$ 
and the line $\lineseg DC$ is parallel to $\lineseg EB$.
(Without the rigidifying gray edge, you would have realizations
which folded half the figure about the line $\lineseg DB$ or $\lineseg EF$.)

\begin{figure}
\BoxedEPSF {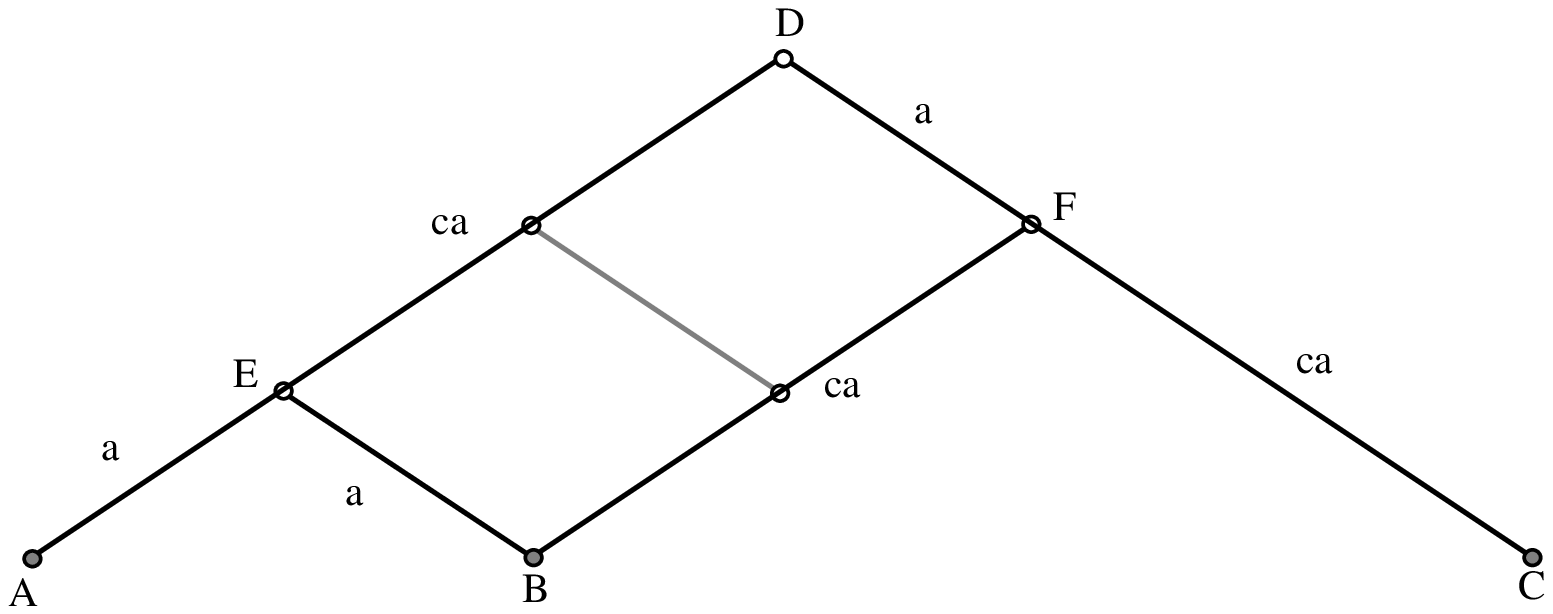 scaled 700}
\capshun 3 { The Pantograph}
\end{figure}

Suppose $c=1$ and the input vertices are $A$ and $C$.
Let the output vertex be $B$.  
Then $\ul$ is \qf\ for $(z,w)\mapsto (z+w)/2$.
Its domain is all $(z,w)$ with $|z-w|\le 4a$ which
can contain any compact set by choosing $a$ big enough.

Next we will take the pantograph and find
\qf\ linkages for $z\mapsto \lambda z$, divided into
three cases: $\lambda >1$, $0<\lambda <1$, and $\lambda <0$.
In all cases the domain is an arbitrarily large ball.

\begin{itemize}
\item
Suppose $B$ is the input vertex, $C$ is 
the output vertex, and $A$ is fixed at 0.
Then the linkage is \qf\ for $z\mapsto (1+c)z$
with domain $|z|\le 2a$.
\item
Suppose $C$ is the input vertex, $B$ is 
the output vertex, and $A$ is fixed at 0.
Then the linkage is \qf\ for $z\mapsto z/(1+c)$
with domain $|z|\le 4a$.
\item
Suppose $A$ is the input vertex, $C$ is 
the output vertex, and $B$ is fixed at 0.
Then the linkage is \qf\ for $z\mapsto -cz$
with domain $|z|\le 2a$.
\end{itemize}

\subsection{Inversion through a circle}

Before constructing the remaining \qfl s,
we will find a \qfl\ for inversion through a circle,
$z\mapsto t^2z/|z|^2$.
This is shown in Figure \ref{fig:4}.
The linkage at the left is the full linkage,
the one at the right just has the basics.
The extra vertices and edges are only needed to
rigidify $BDCE$ and eliminate some degenerate configurations
which occur if $B$ and $C$ coincide.

\begin{figure}
\BoxedEPSF {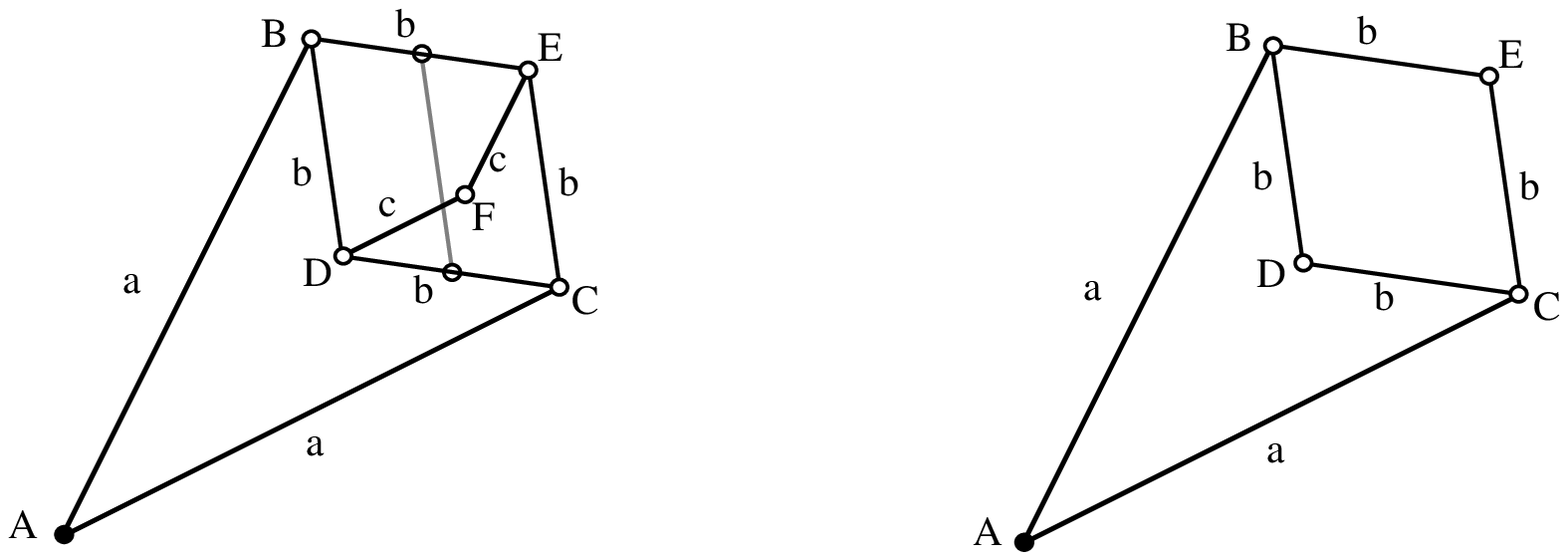 scaled 700}
\capshun 4 { The \PC\ Inversor}
\end{figure}

We fix $A$ at 0, set $t^2=a^2-b^2$, $c<b<a$, let the input vertex be
$D$ and the output be $E$.
Let us see why $\ul$ is \qf\ for $z\mapsto t^2/\overline z=t^2z/|z|^2$.
If $\varphi \in \conf(\ul)$ note that $\varphi (D)$ is a real multiple of $\varphi (E)$.
This follows from the fact that the lines from $\varphi (A)$, $\varphi (D)$, and $\varphi (E)$
to the midpoint of the line segment
$\varphi (B)\varphi (C)$ are all perpendicular to $\varphi (B)\varphi (C)$,
hence $\varphi (A)$, $\varphi (D)$, and $\varphi (E)$ are collinear.
Solving triangles shows that if $s=|\varphi (B)-\varphi (C)|/2$, then 
\begin{eqnarray*}
|\varphi (D)|&=&\sqrt{a^2-s^2}\pm \sqrt{b^2-s^2}\\
 |\varphi (E)|&=&\sqrt{a^2-s^2}\mp \sqrt{b^2-s^2}
 \end{eqnarray*}
from which we see that $|\varphi (D)|\,|\varphi (E)|=(a^2-b^2)$.
So $\ul$ is \qf\ for $z\mapsto t^2z/|z|$.

To see the domain, note that $\sqrt{b^2-c^2}\le s\le b$,
so by the above,  
$$\sqrt{a^2-b^2+c^2}-c\le |\varphi (D)|\le \sqrt{a^2-b^2+c^2}+c$$
So the domain is the annulus between the circles of radius
$\sqrt{t^2+c^2}\pm c$.

\subsection{How to square}

Now let us find a \qfl\ for $z\mapsto z^2$ with domain containing $|z|\le r$.
Note that if $h(z)=t^2z/| z|^2$ then
$$t^2-th((h(t+z)+h(t-z))/2)=z^2$$
Suppose we take a \qfl\ as above for $h$ with $t=4r$ and $c=3r$,
then the domain is $2r\le |z|\le 8r$.
In particular, if $|z|\le r$, then $t+z$, $t-z$, and $(h(t+z)+h(t-z))/2$
are all well within the domain.
So by composition, we get a \qfl\ for $z\mapsto z^2$ 
with domain containing $|z|\le r$.

\subsection{Drawing a straight line}

So finally we are left with finding a \qfl\ for complex conjugation.
Our first step is to find a linkage so that some vertex is constrained to 
lie in a line segment.
This linkage $\ul$ is obtained by taking the input of a quasifunctional linkage for
inversion through a circle and forcing this input to lie in a circle
going through the origin.
But when we invert a circle through the origin, we get a
straight line.
Now it is just a matter of translating and rotating it and rescaling,
to make this line be any interval on the real axis.
This linkage $\ul$ is shown in Figure \ref{fig:5}.

\begin{figure}
\BoxedEPSF {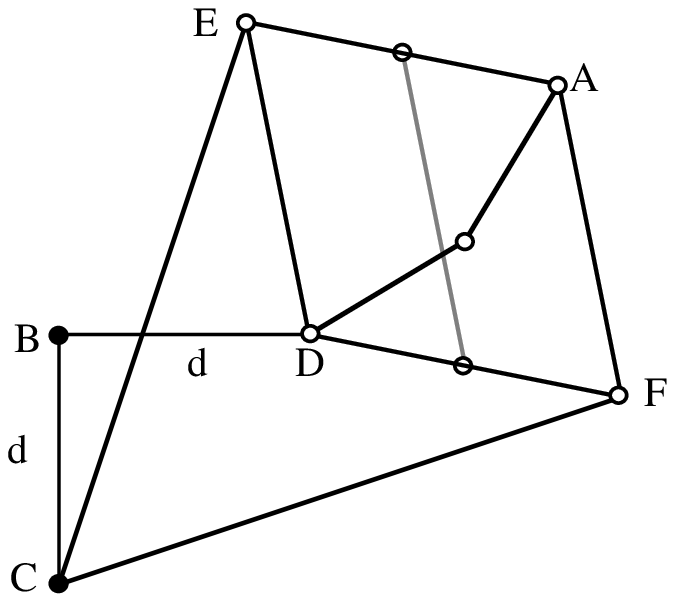 scaled 700}
\capshun {5} { With $B$ and $C$ fixed, $A$ will trace out a straight line}
\end{figure}

We fix $C$ at 0 and fix $B$ at $d\sqrt{-1}$.
Let us first see what $\sconf(\ul,\{D\})$ is.
If $\varphi \in \conf(\ul)$, then we know from the analysis of Figure \ref{fig:4} that
$\varphi (D)$ is in some annulus $0<a\le |z|\le b$.
On the other hand, because of the edge $\lineseg BD$ fixed at $B$,
we must also have $|\varphi (D)-d\sqrt{-1}|=d$.
Thus $\sconf(\ul,\{D\})$ is the intersection of the annulus 
$a\le |z|\le b$ with the circle $|z-d\sqrt{-1}|=d$.
So as long as we choose $d$ so $a<2d<b$, we know this intersection
is an arc of the circle.
But then $\sconf(\ul,\{A\})$ is the inversion of this arc,
which is a straight line segment.

By rescaling all side lengths by a fixed multiple, we may make this line segment
as long as we wish.
By translating and rotating the positions of the fixed points
$C$ and $B$, we may translate and rotate this line segment
to any other line segment with the same length.
So we may construct such a linkage so that $\sconf(\ul,\{A\})$
is any line segment in ${\mathbb C}$ that we wish.

\subsection{Complex conjugation}

We are now ready to construct a quasifunctional linkage for
complex conjugation.
Our first step is to pick two linkages $\ul'$ and $\ul''$ as above
with vertices $A$ and $B$ so that 
$\sconf(\ul',\{A\})=[a,b]$ and $\sconf(\ul'',\{B\})=[-b,-a]$
where $0<a<b$.
We then insert a rigidified square between $A$ and $B$ as shown
in Figure \ref{fig:6}.

\begin{figure}
\BoxedEPSF {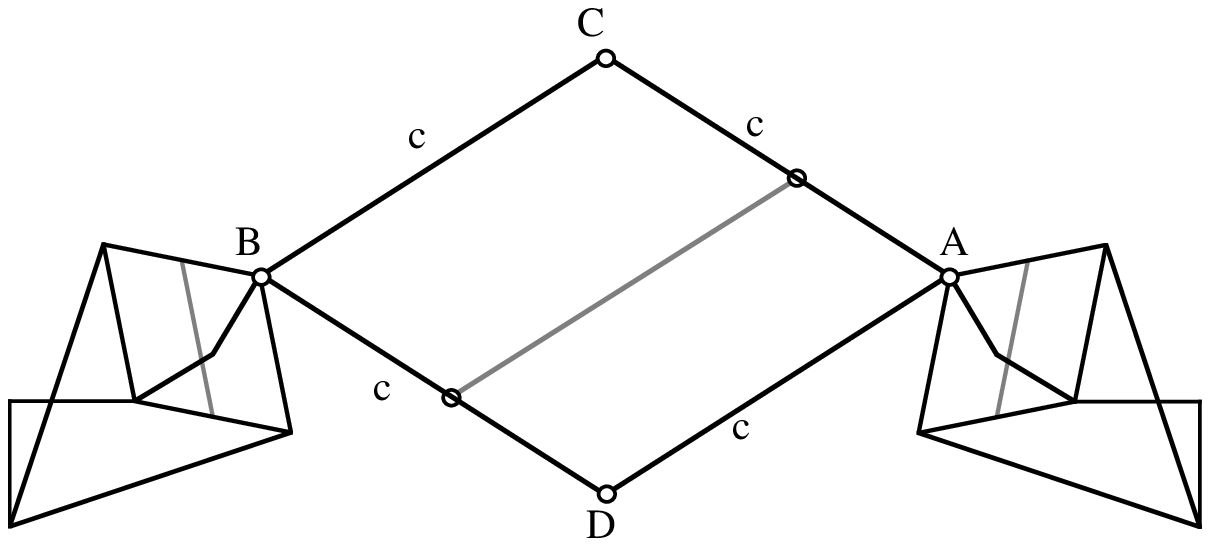 scaled 700}
\capshun {6} { Complex Conjugation}
\end{figure}

Note then that if $C$ is the input vertex and $D$ is the output
vertex, then $\ul$ is \qf\ for $z\mapsto \overline z$.
If we choose $a$, $b$, and $c$ so that $b-r>c>a+r$
then the domain will contain all $z$ with $|z|\le r$.

\section{Proofs of Theorems \ref{thm:a} and \ref{thm:b}}

Now that we have finished proving \thmc,
we may proceed with proving \thma\ and \thmb.
But first we will need a few lemmas.

Recall that a continuous map $f\colon X\to Y$ is proper if 
$f^{-1}(K)$ is compact whenever $K\subset Y$ is compact.
One useful property of proper maps is that if $A\subset X$,
and $X$ and $Y$ are locally compact and Hausdorff,
then $f(\Cl(A))=\Cl(f(A))$. 
Here $\Cl(A)$ stands for the closure of $A$.

I imagine there is a more elementary proof of the
following Lemma which does not use resolution of singularities,
but in any case we have:

\begin{lem}\label{lem:4}
Suppose $X\subset {\mathbb R}^n$ is a semialgebraic set.
Then there is a real algebraic set $Y$ and
a proper polynomial map $q\colon Y\to {\mathbb R}^n$
so that $q(Y)$ is the closure of $X$.
\end{lem}

\begin{proof}
First, we have an algebraic set $X'\subset {\mathbb R}^n\times {\mathbb R}^m$
so that if $p'\colon X'\to {\mathbb R}^n$ is induced by projection,
then $p'(X')=X$.
This is well known and easily illustrated by example.
If $X=\{ x\in {\mathbb R}^n \mid p(x)=0, q(x)\ge 0, r(x)>0 \}$
for polynomials $p$, $q$, and $r$, then just let
$$X'=\{ (x,y,z)\in {\mathbb R}^n\times {\mathbb R}^2 \mid p(x)=0, y^2=q(x), z^2r(x)=1 \}$$

Next, let $X''\subset {\mathbb R}^n\times {\mathbb R}{\mathbb P}^m$ be the Zariski closure of $X'$,
i.e., the smallest real algebraic subset of ${\mathbb R}^n\times {\mathbb R}{\mathbb P}^m$
which contains $X'$.
(${\mathbb R}{\mathbb P}^m$ is  real projective $m$ space.)
Let $p''\colon X''\to {\mathbb R}^n$ be induced by projection.
Note now that $p''$ is proper,
since it is a restriction of the proper map ${\mathbb R}^n\times {\mathbb R}{\mathbb P}^m\to {\mathbb R}^n$
to a closed subset.

It may be true that $X''$ is bigger than the closure
of $X'$, and $p''(X'')$ is bigger than the closure of $X$.
So we must deal with this possibility.
\footnote
{For an example where this occurs, consider the 
case 
\begin{eqnarray*}
X&=&\{ (x,y,z)\in {\mathbb R}^3 \mid x^2<zy^2 \}\\
X'&=&\{ (x,y,z,w)\in {\mathbb R}^4 \mid w^2(zy^2-x^2)=1 \}
\end{eqnarray*}
One can show that 
$X''=\{ (x,y,z,[u:v])\in {\mathbb R}^3\times {\mathbb R}{\mathbb P}^1 \mid 
u^2(zy^2-x^2)=v^2 \}$.
But then points $(0,0,z,[1:0])$ are in $X''$ but not $\Cl(X')$
if $z<0$.}

Let $S$ be the singular set of $X'$ and let $T=X''-X'=X''\cap ({\mathbb R}^n\times {\mathbb R}{\mathbb P}^{m-1})$.
By resolution of singularities (c.f., \cite{H} or \cite{BM}), 
we know that there is a proper map
$p'''\colon X'''\to X''$ so that $X'''$ is nonsingular and $(p''')^{-1} (S\cup T)$
is a union of divisors with normal crossings.
In particular,
$(p''')^{-1} (X'-S)$ is dense in $X'''$, so $p''p'''(X''')$ is the closure
of $p'(X'-S)$.

But by induction on dimension, there is an algebraic set $Y'$
and a proper polynomial map $q'\colon Y'\to {\mathbb R}^n$ so that $q'(Y')$ is the
closure of $p'(S)$.
Letting $Y$ be the disjoint union $X'''\cup Y'$ and letting
$q$ be $p'''\cup q'$, we are done.
\end{proof}

Next we prove an important special case of \thma.

\begin{lem}\label{lem:5}
Suppose $X\subset {\mathbb C}^n$ is a compact semialgebraic set.
Then there is a linkage $\ul$ and a $W\subset \vrtcs(\ul)$
so that $\sconf(\ul,W)=X$.
\end{lem}

\begin{proof}
By Lemma \ref{lem:4}, there is a real proper polynomial map $q\colon Y\to {\mathbb C}^n$
where $Y$ is a real algebraic set and $q(Y)=X$.
By properness of $q$, we know that $Y$ must be compact.
By replacing $Y$ with the graph of $q$,
we may as well suppose that 
$Y\subset {\mathbb C}^n\times {\mathbb C}^m$ for some $m$,
so that $q$ is induced by projection ${\mathbb C}^n\times {\mathbb C}^m\to {\mathbb C}^n$.
Let $p\colon {\mathbb C}^n\times {\mathbb C}^m\to {\mathbb C}$ be a real polynomial so that $Y=p^{-1}(0)$.

By \thmc, there is a  \qfl\ $\ul'$
for $p$ 
with distinct input and output vertices and
whose domain contains $Y$.
Construct a linkage $\ul$ by taking  $\ul'$ 
and fixing its output vertex to 0.
Let $W$ be the set of the first $n$ input vertices
and let $U$ be the set of all input vertices.
Since the output vertex of $\ul'$ is fixed to
0 and $\ul'$ is \qf\ for $p$ 
we know that the input vertices must always lie on $Y=p^{-1}(0)$.
So $\sconf(\ul,U)=Y$.
Consequently,
$\sconf(\ul,W)$ is the projection of $Y$ to
${\mathbb C}^n$ which is $X$.
\end{proof}

\begin{lem}\label{lem:6}
Suppose $\ul$ is the disjoint union of two
linkages $\ul'$ and $\ul''$,  and $W\subset \vrtcs(\ul)$.
Then after reordering $W$ so that
the vertices in $W\cap \vrtcs(\ul')$ come first, we have
$$\sconf(\ul,W)=\sconf(\ul',W\cap \vrtcs(\ul'))\times
 \sconf(\ul'',W\cap \vrtcs(\ul''))$$
\end{lem}

\begin{proof} Any planar realization of $\ul$ 
is a pair of realizations of $\ul'$ and $\ul''$,
and vice versa.
\end{proof}

\begin{lem}\label{lem:7}
Suppose $Z\subset {\mathbb C}^k$ is invariant under the diagonal action of $\Euc(2)$.
Let $Z_0=Z\cap ({\mathbb C}^{k-1}\times 0)$ be the points in $Z$ with last coordinate 0.
Then:
\begin{enumerate}
\item $Z$ is virtually compact if and only if $Z_0$ is compact.
\item $Z_0$ is invariant under the diagonal action of $O(2)$.
\item  $Z_0/O(2)$ is homeomorphic to $Z/\Euc(2)$.
\item If $Y\subset {\mathbb C}^k$ is invariant under the diagonal action of $\Euc(2)$
and $Z_0=Y\cap ({\mathbb C}^{k-1}\times 0)$, then $Y=Z$.
\end{enumerate}
\end{lem}

\begin{proof}
If $z\in Z_0$ and $\beta \in O(2)$ then $\beta (Z)\in Z$ and $\beta (z)\in {\mathbb C}^{k-1}\times 0$,
so $\beta (z)\in Z_0$, and 2 is shown.

We have a map $f\colon Z\to Z_0$ given by $f(z_1,\ldots ,z_k)=(z_1-z_k,\ldots ,z_{k-1}-z_k,0)$.
Then $f$ induces a continuous map $f'\colon Z/\Euc(2)\to Z_0/O(2)$.
Likewise the inclusion $Z_0\subset Z$ induces a continuous map $g'\colon Z_0/O(2)\to Z/\Euc(2)$
and these maps are inverses of each other.
So 3 is proven.

To see 4, note that if $\Tran(2)\subset \Euc(2)$ is the subgroup of
translations, then both $Z$ and $Y$ are the union of $\Tran(2)$ orbits of points of $Z_0$.
So $Y=Z$.

Now suppose that $Z$ is virtually compact.
Then by 3 we know that $Z_0/O(2)$ is compact.
Suppose $z\not\in Z_0$ but $z$ is in the closure of $Z_0$.
Let $K$ be the $O(2)$ orbit of $z$.
For any $\delta >0$ let $U_\delta $ be the set of points in ${\mathbb C}^k$ with distance
greater than $\delta $ from $K$.
Each $U_\delta $ is $O(2)$ invariant, so we get an open cover
$\{ (U_\delta \cap Z_0)/O(2) \}$ of $Z_0/O(2)$.
By compactness of $Z_0/O(2)$ we know that $U_\delta \supset Z_0$ for some $\delta >0$,
contradicting $z$ being in the closure of $Z_0$.
So $Z_0$ is closed.
For any $r>0$ let $B_r$ be the open ball of radius $r$ around 0.
We get an open cover $\{ (B_r\cap Z_0)/O(2) \}$ of $Z_0/O(2)$.
So $B_r\supset Z_0$ for some $r$ and thus $Z_0$ is bounded.
So $Z_0$ is compact.

On the other hand, if $Z_0$ is compact, then $Z_0/O(2)$ is compact,
so $Z$ is virtually compact, so 1 is shown.
\end{proof}

We are now ready to prove \thmb.

\begin{proof} (of \thmb)
We first prove 5-8.
So suppose $\ul$ is a connected linkage with $m$ fixed vertices
and $W\subset \vrtcs(\ul)$.
By the Tarski-Seidenberg theorem \cite{S},
we know $\sconf(\ul,W)$ is a semialgebraic set
since it is a projection of the algebraic set $\conf(\ul)$.

Note that if $\beta \in \Euc(2)$ and $\beta $ fixes the images
of all fixed vertices of $\ul$, and $\varphi \in \conf(\ul)$, 
then $\beta \varphi \in \conf(\ul)$ also.
Consequently $\beta (\sconf(\ul,W))\subset \sconf(\ul,W)$.
So if $m=0$,  $\sconf(\ul,W)$ is invariant under $\Euc(2)$.
If $m=1$, $\sconf(\ul,W)$ is invariant under the
subgroup $G$ of $\Euc(2)$ which fixes the image of the
fixed vertex of $\ul$.
If $m=2$, $\sconf(\ul,W)$ is invariant under the
subgroup $G$ of $\Euc(2)$ which fixes a line
through the images of the
two fixed vertices of $\ul$.
So we have shown everything but compactness.

If $m>0$, let $z_0$ be the image of some fixed vertex
and let $d$ be the sum of the lengths of all
edges of $\ul$.
Then for each $\varphi \in \conf(\ul)$
and each $v\in \vrtcs(\ul)$ we know that 
$\varphi (v)$ lies in the ball of radius $d$ around $z_0$.
So $\conf(\ul)$ is bounded.
But it is also closed since it is an algebraic set.
So $\conf(\ul)$ is compact.
Hence $\sconf(\ul,W)$ is compact since it is the
continuous image of $\conf(\ul)$.

If $m=0$, let $W=\{w_1,\ldots ,w_k\}$.
If $k=0$ then $\sconf(\ul,W)={\mathbb C}^0$ is compact,
so assume that $k>0$.
Let $\ul'$ be the linkage obtained from $\ul$
by fixing the vertex $w_k$ at 0.
Let $Z_0=\sconf(\ul',W)$.
Note that $Z_0$ is compact and $O(2)$ invariant by 6.
Also $Z_0=\sconf(\ul,W)\cap {\mathbb C}^{k-1}\times 0$.
So by Lemma \ref{lem:7}-1, $Z$ is virtually compact so  5 is proven.
Note that $\sconf(\ul,W)$ must be noncompact
since it is invariant under $\Tran(2)$.

Now we will prove 2, 3, and 4.
Note in these cases that $Z$ is compact.
After replacing $Z$ by $\beta (Z)$ for some $\beta \in \Euc(2)$,
 we may assume in case 2 that $G=O(2)$, 
 and may assume in case 3 that
$Z$ is invariant under complex conjugation.
By Lemma \ref{lem:5},  we may find a linkage $\ul'$ and a $W\subset \vrtcs(\ul')$ so that
$\sconf(\ul',W)=Z$.
Throw away all connected components of $\ul'$ which do not contain
any vertices in $W$ or any fixed vertices, doing so does not change $\sconf(\ul',W)$.
By adding some isolated fixed vertices to $\ul'$
if necessary, we may assume that there is a vertex fixed at 0,
another fixed at 1,  a third fixed at $\sqrt{-1}$, and a fourth fixed at $-1-\sqrt{-1}$.
Adding an isolated fixed vertex to $\ul'$ does not change
$\sconf(\ul',W)$.

Let the fixed vertices of $\ul'$ be $\{v_1,\ldots ,v_k\}$
where $v_i$ is fixed to the point $z_i$.
We may suppose $z_1=0$, $z_2=1$,  $z_3=\sqrt{-1}$, and $z_4=-1-\sqrt{-1}$.
For each pair $i,j$ with $z_i\neq z_j$ put in an edge 
$\lineseg{v_i}{v_j}$ of length $|z_i-z_j|$,
if it is not already there.
This will not change $\sconf(\ul',W')$.
%
%
%
%
%
Note we did not attempt to add any zero length edges,
which would not be allowed.

Let $\ul''$ be obtained from $\ul'$ by only fixing the vertices
$v_1$, $v_2$, and $v_3$.
We claim that $\sconf(\ul',W)=\sconf(\ul'',W)$.
One inclusion $\sconf(\ul',W)\subset \sconf(\ul'',W)$ is trivial.
So let us see the other inclusion.
Pick any $\varphi \in \conf(\ul'')$. 
We claim that in fact $\varphi (v_i)=z_i$ for all $i$.
To see this, note first that two different points in ${\mathbb C}$ can't
have the same distances from three noncollinear points.
Consequently $\varphi '(v_4)=z_4$ since the three edges 
$\lineseg{v_i}{v_4}$, $i=1,2,3$
have lengths $|z_4-z_i|$, so $|\varphi '(v_4)-z_i|=|\varphi '(v_4)-\varphi '(v_i)|=|z_4-z_i|$.
Then for any $j>4$, there are edges in $\ul$ from $v_j$ to at least three
of the $v_i$, $i\le 4$, and consequently $\varphi '(v_j)=z_j$ since any three
of the $z_i$, $i\le 4$ are noncollinear.
Consequently, $\varphi \in \conf(\ul')$.
So $\varphi |_{W}\in \sconf(\ul',W)$, and we have shown that
$\sconf(\ul',W)=\sconf(\ul'',W)$.

We claim that $\ul''$ is connected.
All fixed vertices of $\ul$ are in the same
connected component since they are all connected
by edges. 
We also threw out any components without points of $W$.
So any other components must have no fixed vertices and
must contain points of $W$.
We saw from the proof of 5 above that semiconfiguration spaces
of such components are noncompact.
Hence by Lemma \ref{lem:6}, the configuration space of $\ul''$
would be noncompact, but it is not.
So $\ul''$ is connected and so 4 is proven.

Let us now prove 3.
Let $\ul$ be obtained from $\ul''$ by only fixing the vertices $v_1$ and $v_2$,
and not fixing $v_3$.
We claim that $\sconf(\ul'',W)=\sconf(\ul,W)$.
Again, one inclusion $\sconf(\ul'',W)\subset \sconf(\ul,W)$ is trivial.
So let us see the other inclusion.
Pick any $\varphi \in \conf(\ul)$. 
Now $|\varphi (v_3)|=|z_3|=1$, and $|\varphi (v_2)-\varphi (v_3)|=|z_2-z_3|=\sqrt 2$,
so the triangles $z_1z_2z_3$ and $z_1z_2\varphi (v_3)$ are congruent,
and hence $\varphi (v_3)=\pm \sqrt{-1}$.
If $\varphi (v_3)=\sqrt{-1}$ then $\varphi \in \conf(\ul'')$ 
so $\varphi |_W\in \sconf(\ul'',W)$.
If $\varphi (v_3)=-\sqrt{-1}$ then the complex conjugate 
$\overline \varphi \in \conf(\ul'')$ so 
$\overline \varphi |_W\in \sconf(\ul'',W)=Z$,
but then $\varphi |_W\in \sconf(\ul'',W)$ since $Z$ is invariant under
complex conjugation.

Now let us prove 2.
Let $\ul$ be obtained from $\ul'$ by only fixing the vertex $v_1$ to 0,
and not fixing any of the other vertices of $\ul'$.
We claim that $\sconf(\ul'',W)=\sconf(\ul,W)$.
One inclusion $\sconf(\ul'',W)\subset \sconf(\ul,W)$ is trivial.
So let us see the other inclusion.
Pick any $\varphi \in \conf(\ul)$. 
Now $|\varphi (v_2)|=|z_2|=1$, $|\varphi (v_3)|=|z_3|=1$, and $|\varphi (v_2)-\varphi (v_3)|=|z_2-z_3|=\sqrt 2$,
 so there is a  $\beta \in O(2)$ so that
$\beta (\varphi (v_2))=z_2=1$ and $\beta (\varphi (v_3))=z_3=\sqrt{-1}$.
For convenience, let $\varphi '=\beta \compose \varphi $.
Note that $\varphi '\in \conf(\ul'')$.
So $\varphi '|_{W}\in \sconf(\ul'',W)=Z$.
By $O(2)$ invariance of $Z$, we know that $\beta ^{-1}\compose \varphi '|_{W}\in Z$ also.
But $\beta ^{-1}\compose \varphi '|_{W}=\varphi |_{W}$, so $\varphi |_{W}\in Z$.
So we have shown that $\sconf(\ul,W)\subset \sconf(\ul'',W)$,
and hence $\sconf(\ul,W)=Z$.

Finally, let us prove 1.
Let $Z_0=Z\cap {\mathbb C}^{n-1}\times 0$.
By Lemma \ref{lem:7}-1 and 2, $Z_0$ is compact and $O(2)$ invariant.
Hence by 2 there is a $\ul'$ and $W$ so that 
$Z_0=\sconf(\ul',W)$ and $\ul'$ has only one fixed vertex,
fixed at 0.
Form $\ul$ from $\ul'$ by unfixing this vertex.
Then $\sconf(\ul,W)\cap {\mathbb C}^{n-1}\times 0=\sconf(\ul',W)=Z_0$.
So $\sconf(\ul,W)=Z$ by Lemma \ref{lem:7}-4.
\end{proof}

We can now prove \thma.

\begin{proof} (of \thma)
By Lemma \ref{lem:6}, it suffices to prove the equivalence of 2(a) and 2(b).
The equivalence of 1(a) and 1(b) will follow by applying the equivalence
of 2(a) and 2(b) to each connected component of $\ul$ which
contains vertices of $W$.

We know from \thmb-8 and 5 that 2(a) implies 2(b),
and along with Lemma \ref{lem:6} we get all the extra bits of information.

So it remains to prove 2(b) implies 2(a).
If $X$ is compact, this follows from \thmb -4.
Suppose now $X$ is not compact,
so $X$ is invariant under the action of $\Eu(2)$ 
with compact quotient.
Then 2(a) follows from \thmb-1.
\end{proof}


\begin{thebibliography}{99}

\bibitem{AK} S.~Akbulut and H.~King, \textit{Topology of real algebraic sets},
MSRI Publ.~25, Springer-Verlag (1992).


\bibitem{BM} E.~Bierstone and P.~Milman, \textit{Canonical desingularization 
in characteristic zero by blowing up the maximum strata
of a local invariant}, Inventiones Math.~128 (1997), 207-302.

\bibitem{H} H.~Hironaka,  \textit{Resolution of singularities of an algebraic variety
over a field of characteristic zero},
Annals of Math.~79 (1964), 109-326.

\bibitem{K} H.~King, \textit{Planar linkages and algebraic sets},
 preprint, Math.AG/9807023.

\bibitem{KM} M.~Kapovich and J.~Millson, \textit{Universality theorems for 
configuration spaces of planar linkages}, preprint.

\bibitem{S} A.~Seidenberg, \textit{A new decision method for elementary algebra},
Annals of Math.~60 (1954), 365-374.

\end{thebibliography}
\end{document}